\documentclass[11pt,amsfonts]{article}
\usepackage{graphicx}
\usepackage{latexsym}
\usepackage{amssymb}
\usepackage{amsmath}
\usepackage{enumerate}
\usepackage{layout}
\usepackage{eufrak}
\usepackage{color}
\newtheorem{prop}{Proposition}
\newtheorem{lemma}{Lemma}

\newtheorem{theorem}{Theorem}
\newtheorem{remark}{Remark}

\def\real{{\mathord{{\rm I\kern-2.8pt R}}}}        
\def\inte{{\mathord{{\rm I\kern-2.8pt N}}}}

\def\sZZ{{\rm Z\kern-2.8ptem{}Z}}

\def\z{{\mathchoice
  {\sZZ}
  {\sZZ}
  {\rm Z\kern-0.30em{}Z}
  {\rm Z\kern-0.25em{}Z} }}
\def\sQQ{{\kern 0.27em \vrule height1.45ex width0.03em depth0em
          \kern-0.30em \rm Q}}
\def\qu{{\mathchoice
    {\sQQ}
    {\sQQ}
  {\kern 0.225em \vrule height1.05ex width0.025em depth0em \kern-0.25em \rm Q}
  {\kern 0.180em \vrule height0.78ex width0.020em depth0em \kern-0.20em \rm Q}
        }}
\def\sCC{{\kern 0.27em \vrule height1.45ex width0.03em depth0em
          \kern-0.30em \rm C}}
\def\complex{{\mathchoice
    {\sCC}
    {\sCC}
  {\kern 0.225em \vrule height1.05ex width0.025em depth0em \kern-0.25em \rm C}
  {\kern 0.180em \vrule height0.78ex width0.020em depth0em \kern-0.20em \rm C}
        }}

\newcommand{\ba}{\begin{array}}
\newcommand{\ea}{\end{array}}
\newcommand{\be}{\begin{equation}}
\newcommand{\ee}{\end{equation}}
\newcommand{\bea}{\begin{eqnarray}}
\newcommand{\eea}{\end{eqnarray}}
\newcommand{\beaa}{\begin{eqnarray*}}
\newcommand{\eeaa}{\end{eqnarray*}}

\def\z{\zeta}

\font\tenmath=msbm10 \font\sevenmath=msbm7 \font\fivemath=msbm5
\newfam\mathfam \textfont\mathfam=\tenmath
\scriptfont\mathfam=\sevenmath \scriptscriptfont\mathfam=\fivemath

\def \={{\buildrel {\rm (law)} \over =}}

\def\qed{ \hfill \vrule width.25cm height.25cm depth0cm\smallskip}

\newcommand{\basa}{\begin{assumption}}
\newcommand{\easa}{\end{assumption}}

\newcommand{\bas}{\begin{assum}}
\newcommand{\eas}{\end{assum}}

\def\limsup{\mathop{\overline{\rm lim}}}
\def\liminf{\mathop{\underline{\rm lim}}}

\newcommand{\ignore}[1]{}
\begin{document}
\allowdisplaybreaks
\renewcommand{\thefootnote}{\fnsymbol{footnote}}
\renewcommand{\thefootnote}{\fnsymbol{footnote}}
\title{Stein's method for invariant measures of diffusions via Malliavin calculus \footnote{Dedicated to the memory of Constantin Tudor. }}
\author{ Seiichiro Kusuoka $^{1}$
\footnote{ Research Fellow of the Japan Society for the Promotion of Science}$\qquad $%
Ciprian A. Tudor $^{2,}$
\footnote{ Partially supported by the ANR grant "Masterie" BLAN 012103. Associate member of the team Samos, Universit\'e de Panth\'eon-Sorbonne Paris 1.  This work has been partially completed
 when the second author has visited Keio University. He acknowledges generous support from Japan Society for the Promotion of Science. }\vspace*{0.1in} \\
$^{1}$Department of Mathematics, Kyoto University \\
Kyoto 606-8502, Japan  \\
kusuoka@math.kyoto-u.ac.jp \vspace*{0.1in} \\
 $^{2}$ Laboratoire Paul Painlev\'e, Universit\'e de Lille 1\\
 F-59655 Villeneuve d'Ascq, France.\\
 \quad tudor@math.univ-lille1.fr\vspace*{0.1in}}
\maketitle

\begin{abstract}
\noindent Given a random variable $F$ regular enough in the sense of the Malliavin calculus, we are able to measure the distance between its law and any probability measure with a density function which is continuous, bounded, strictly positive on an interval in the real line and admits finite variance. The bounds are given in terms of the Malliavin derivative of $F$. Our approach is based  on the theory of It\^o diffusions and the stochastic calculus of variations. Several examples are considered in order to illustrate our general results.
\end{abstract}

\vskip0.2cm

{\bf 2010 AMS Classification Numbers:}  60F05, 60H05, 91G70.

 \vskip0.2cm


\section{Introduction}
Let $(\Omega, {\mathcal{F}}, P)$ be a probability space and $(W_{t})_{t \geq 0}$ a Brownian motion on this space. Let $F$ be a random variable defined on $\Omega$ which is differentiable in the sense of the Malliavin calculus.
Then using  the so-called Stein's method introduced by Nourdin and Peccati in \cite{NoPe1} (see also \cite{NoPe2} and \cite{NoPe3}), it is possible to measure the distance between the law of $F$ and the standard normal law $N(0,1)$. This distance can be defined in several ways, e.g. the Kolmogorov distance, the Wasserstein distance, the total variation distance or the Fortet-Mourier distance. More precisely we have, if ${\mathcal{L}}(F)$  denotes the law of $F$,
\begin{equation}\label{b1}
d({\mathcal{L}}(F),N(0,1)) \leq c\sqrt{ \mathbf{E}\left( 1-\langle DF , D(-L) ^{-1} F \rangle _{L^{2}([0,1])}\right)^{2}}.
\end{equation}
Here $D$ denotes the Malliavin derivative with respect to $W$, and $L$ is the generator of the Ornstein-Uhlenbeck semigroup. We will explain in the next section  how these operators are defined. The constant $c$  in (\ref{b1}) depends on the distance considered: is equal to 1 in the case of the Kolmogorov distance as well as in the case of the Wasserstein distance, $c=2$ for the  total variation distance and $c=4$ in the case of the Fortet-Mourier distance.

The Stein's method related with Malliavin calculus has been generalized to Gamma and Pearson distributions (see \cite{NoPe4} and \cite{EV}). Our purpose is to extend the above bound to a more general class of probability distributions. Actually, given a random variable $F$ regular enough in the sense of the Malliavin calculus, we are able to measure the distance between its law and any probability measure with a density function which is continuous, bounded, strictly positive on an interval in the real line and admits finite variance. Our approach is based on the following idea. Consider a given probability density function $p$ which is continuous, bounded, with finite variance, it is strictly positive on an interval $(l,u)$ ($-\infty \leq l<u\leq \infty)$ and it vanishes outside the interval $(l,u)$. This class includes almost all the commonly used probability densities. Then we will construct an ergodic It\^o diffusion which admits $p$ as invariant density. The procedure to construct such a diffusion is explained in \cite{BSS} (see also \cite{AS} or \cite{Reinert}) and it will be recalled in Section \ref{generaltheory} in our paper. The diffusion associated with the density $p$ (by associated with $p$, we mean that it admits $p$ as invariant  density) has the form
\begin{equation*}
dX_{t} = -(X_{t}-m) dt + \sqrt{a(X_{t})} dW_{t}
\end{equation*}
where $m$  is the mean of $p$, $(W_{t})_{t\geq 0}$ is a standard Wiener process and the diffusion coefficient $a$ can be written explicitly in terms of the function $p$ and of the constant $m$ (see formula (\ref{difcoef}) in Section \ref{generaltheory}).
Then we will consider the generator of the diffusion $X$ above and we will use the properties of this operator and the integration by parts on Wiener space in the spirit of \cite{NoPe1} in order to obtain a bound between the law of an arbitrary random variable $F$ and the law with density $p$. This bound, that we will call Stein's bound, will involve the Malliavin derivative of $F$ and it reduces to the results contained in \cite{NoPe1} and \cite{NoPe4} in the particular cases of the Gaussian and Gamma laws.

We mention that there already exists a scientific literature relating Stein's method and Berry-Ess\'een bound to diffusion theory. We refer to the survey \cite{Reinert} and the reference therein (see among others \cite{Barbour}, \cite{Ste}) or to the recent reference \cite{EV}. The novelty of our approach is related to the use of the Malliavin derivatives. This approach has several advantages. Firstly, it is unitary, in the sense the Stein's bound we give is available for several metrics between the laws of random variables. Secondly, it is rather complete since, as we explain in Section 3, it applies to almost every absolutely continuous probability distribution. And thirdly, it allows us to characterize the random variables with a given probability distribution in terms of its Malliavin derivatives (see Theorem  \ref{char}).

Our paper is structured in the following way: Section 2 contains the basic notion on Malliavin calculus. In Section 3 we construct our general theory to derive the Berry-Ess\'een bound for the distance between an arbitrary random variable and a given probability measure. In Sections 4 and we consider several examples (the uniform distribution, the Pareto distribution, the Laplace distribution etc) used to illustrate our bound.

\section{Preliminaries}
 This paragraph is devoted to  introduce the elements from stochastic analysis that will be used  in the paper. Consider ${\mathcal{H}}$ a real separable Hilbert space and $(B (\varphi), \varphi\in{\mathcal{H}})$ an isonormal Gaussian process on a probability space $(\Omega, {\cal{A}}, P)$, which is a centered Gaussian family of random variables such that $\mathbf{E}\left( B(\varphi) B(\psi) \right)  = \langle\varphi, \psi\rangle_{{\mathcal{H}}}$. Denote  $I_{n}$ the multiple stochastic integral with respect to
$B$ (see \cite{N}). This $I_{n}$ is actually an isometry between the Hilbert space ${\mathcal{H}}^{\odot n}$(symmetric tensor product) equipped with the scaled norm $\frac{1}{\sqrt{n!}}\Vert\cdot\Vert_{{\mathcal{H}}^{\otimes n}}$ and the Wiener chaos of order $n$ which is defined as the closed linear span of the random variables $H_{n}(B(\varphi))$ where $\varphi\in{\mathcal{H}}, \Vert\varphi\Vert_{{\mathcal{H}}}=1$ and $H_{n}$ is the Hermite polynomial of degree $n\geq 1$ given by
\begin{equation*}
H_{n}(x)=\frac{(-1)^{n}}{n!} \exp \left( \frac{x^{2}}{2} \right)
\frac{d^{n}}{dx^{n}}\left( \exp \left( -\frac{x^{2}}{2}\right)
\right), \hskip0.5cm x\in \mathbb{R}.
\end{equation*}
The isometry of multiple integrals can be written as: for $m,n$ positive integers,
\begin{eqnarray}
\mathbf{E}\left(I_{n}(f) I_{m}(g) \right) &=& n! \langle \tilde{f},\tilde{g}\rangle _{{\mathcal{H}}^{\otimes n}}\quad \mbox{if } m=n,\nonumber \\
\mathbf{E}\left(I_{n}(f) I_{m}(g) \right) &= & 0\quad \mbox{if } m\not=n.\label{iso}
\end{eqnarray}
It also holds that $I_{n}(f) = I_{n}\big( \tilde{f}\big)$ where $\tilde{f} $ denotes the symmetrization of $f$ defined by $\tilde{f}%
(x_{1}, \ldots , x_{n}) =\frac{1}{n!} \sum_{\sigma \in {\cal S}_{n}}
f(x_{\sigma (1) }, \ldots , x_{\sigma (n) } ) $.
\\\\
We recall that any square integrable random variable which is measurable with respect to the $\sigma$-algebra generated by $B$ can be expanded into an orthogonal sum of multiple stochastic integrals
\begin{equation}
\label{sum1} F=\sum_{n\geq0}I_{n}(f_{n})
\end{equation}
where $f_{n}\in{\mathcal{H}}^{\odot n}$ are (uniquely determined)
symmetric functions and $I_{0}(f_{0})=\mathbf{E}\left[  F\right]$.
\\\\
Let $L$ be the Ornstein-Uhlenbeck operator
\begin{equation*}
LF=-\sum_{n\geq 0} nI_{n}(f_{n})
\end{equation*}
if $F$ is given by (\ref{sum1}) and it is such that $\sum_{n=1}^{\infty} n^{2}n! \Vert f_{n} \Vert ^{2} _{{\cal{H}}^{\otimes n}}<\infty$.
\\\\
For $p>1$ and $\alpha \in \mathbb{R}$ we introduce the Sobolev-Watanabe space $\mathbb{D}^{\alpha ,p }$  as the closure of
the set of polynomial random variables with respect to the norm
\begin{equation*}
\Vert F\Vert _{\alpha , p} =\Vert (I -L)^{\frac{\alpha }{2}} F \Vert_{L^{p} (\Omega )}
\end{equation*}
where $I$ represents the identity. We denote by $D$  the Malliavin  derivative operator that acts on smooth functions of the form $F=g(B(\varphi _{1}), \ldots , B(\varphi_{n}))$ ($g$ is a smooth function with compact support and $\varphi_{i} \in {{\cal{H}}}$)
\begin{equation*}
DF=\sum_{i=1}^{n}\frac{\partial g}{\partial x_{i}}(B(\varphi _{1}), \ldots , B(\varphi_{n}))\varphi_{i}.
\end{equation*}
The operator $D$ is continuous from $\mathbb{D}^{\alpha , p} $ into $\mathbb{D} ^{\alpha -1, p} \left( {\cal{H}}\right).$
The adjoint of
$D$ is denoted by $\delta $ and is called the divergence (or
Skorohod) integral. It is a continuous operator from $\mathbb{D} ^{\alpha -1, p } \left( {\cal{H}}\right)$ into $\mathbb{D} ^{\alpha , p}$.
The following key relation provides a connection between $D, L$ and $\delta $ and plays an important role in the so-called Stein's method: for any centered random variables $F \in \mathbb{D} ^{1,2}$ it holds that
 $F = LL^{-1}F = -\delta DL^{-1}F$
\section{The general theory}\label{generaltheory}
\subsection{It\^o diffusion with given invariant measure}

In this paragraph we will describe the construction of a diffusion process with given invariant measure $\mu$ that admits a density $p$ with respect to the Lebesque measure. We refer to \cite{BSS} and \cite{AS} for more details and proofs.  Assume that the density $p$ satisfies the following conditions: it is continuous, bounded, admits finite variance and $p$ is strictly positive on the interval $(l,u) $ ($-\infty \leq l<u\leq \infty$) and it is zero outside $(l,u)$. Denote by $m$ the expectation of $\mu$ and consider the stochastic differential equation
\begin{equation}
\label{dif1}
dX_{t} = -( X_{t}-m) dt + \sqrt{ a(X_{t}) } dW_{t}, \hskip0.5 cm t\geq 0
\end{equation}
where $(W_{t}) _{t\geq 0}$ is a standard Brownian motion and the diffusion coefficient is defined by
\begin{equation}\label{difcoef}
a(x)= \frac{2 \int_{l}^{x} ( m-y) p(y) dy}{p(x) } = \frac{ 2 m F(x) -2 \int_{l}^{x} yp(y) dy}{p(x) }, \hskip0.5cm x\in (l,u)
\end{equation}
where $F(x)= \int_{-\infty} ^{x} p(y)dy, x\in \mathbb{R}$  is the distribution function associated with the density $p$. Then the following holds (see Theorem  2.3 in \cite{BSS}):
\begin{description}
\item{$\bullet $ } The stochastic differential equation (\ref{dif1}) with diffusion coefficient given by (\ref{difcoef}) has a unique Markovian weak solution.
   \item{$\bullet $ } The diffusion coefficient $a$ (\ref{difcoef})  is strictly positive for $x\in (l,u)$ and satisfies
   $$\mathbf{E} a(X)= \int_{l}^{u} a(x) p(x) dx <\infty $$
   where $X \sim \mu$ (this notation means that the random variable $X$ follows the law $\mu$).
   \item{$\bullet $ } The solution $X$ to (\ref{dif1}) is ergodic with invariant density $p$.
   \item{$\bullet $ }  If $ -\infty <l $ or $u<\infty$ then (\ref{dif1}) is the only ergodic diffusion with drift $-(x-m)$ and invariant density $p$. If the state space is the whole real line, then (\ref{dif1}) is the only ergodic diffusion with drift $-(x-m)$ and invariant density $p$ such that $\int_{l}^{u} a(x) p(x) dx <\infty$.
\end{description}
Table 1 in \cite{BSS} provided many examples of diffusion associated with a given density $p$. We will  use some of them in our paper (the normal distribution, the Gamma distribution, the uniform distribution, the Beta distribution, the log-normal distribution, the Laplace  distribution, the log-normal distribution) and we will recall the diffusion coefficients associated with this law. Besides these examples, many others can be found in \cite{BSS} and the list is not exhaustive. In principle for any density that satisfies the rather general assumption described at the beginning of this section, one can associate a diffusion process. For some classes of distributions it is not possible to determine an explicit expression for the squared diffusion coefficient $a$ similar to (\ref{difcoef}). In this case, some approximation method can be applied (see Section 3 in \cite{BSS}).
\begin{remark}\label{remark1}
The construction of the diffusion process presented above is a particular case of a more general result. That is, given a density $p$ as above and given a drift coefficient $b$ such that there exists a real number  $k \in (l,u)$ such that $b(x)>0$  for $x \in (l,k)$ and $b(x)<0$ for $x\in (k,u)$, $bp$ is continuous and bounded on $(l,u)$  and
$$\int_{l}^{u} b(x)p(x)dx =0$$
then  the stochastic differential equation
\begin{equation*}
dX_{t}= b(X_{t}) dt + \sqrt{a(X_{t}) } dW_{t}, \hskip0.5cm t\geq 0
\end{equation*}
with $a(x)= \frac{2 \int_{l}^{x} b(y) p(y) dy}{p(x)}$ has a unique Markovian weak solution which is ergodic with invariant density $p$.
\end{remark}
This way is not the only one to construct a diffusion process with a given invariant density. Another method for constructing such a diffusion is given in \cite{BS}.

\subsection{Stein's method of invariant measures associated to one-dimensional second order differential operators}
The purpose of this paragraph is to derive the bounds for the distance between the law of an arbitrary random variable $Y$ and a given continuous probability distribution. We will situate ourselves in the general context described in Remark \ref{remark1}.

Let $S$ be the interval $(l,u) $ ($-\infty \leq l<u\leq \infty$) and $\mu$ be a probability measure on $S$ with a density function $p$ which is continuous, bounded, strictly positive on $S$, and admits finite variance.
Consider a continuous function $b$ on $S$ such that there exists $k \in (l,u)$ such that $b(x)>0$  for $x \in (l,k)$ and $b(x)<0$ for $x\in (k,u)$, $bp$ is bounded on $S$ and
\begin{equation}\label{cond-b}
\int_{l}^{u} b(x)p(x)dx =0.
\end{equation}
Define
\begin{equation}\label{def-a}
a(x):= \frac{2 \int_{l}^{x} b(y) p(y) dy}{p(x)},\quad x\in S.
\end{equation}
Then, the stochastic differential equation:
\begin{equation*}
dX_{t}= b(X_{t}) dt + \sqrt{a(X_{t}) } dW_{t}, \hskip0.5cm t\geq 0
\end{equation*}
has a unique Markovian weak solution,  ergodic with invariant density $p$ (see Remark \ref{remark1}).

In the next result we express the density $p$ in terms of the coefficients $a$ and $b$.
The equation in the following proposition is same as the formula (6.22) on page 241 in \cite{KaTa}.

\begin{prop}\label{prop-ap}Take $c\in S$, then
\[
p(x)= \frac{p(c)a(c)}{a(x)}e^{\int _c^x \frac{2b(y)}{a(y)}dy}, \quad x\in S.
\]
\end{prop}
{\bf Proof: }
Relation (\ref{def-a}) implies that
\[
a(x)p(x) = 2\int _l^x b(y)p(y)dy.
\]
Note that the left-hand side of this equality is differentiable, since the right-hand side is differentiable.
Differentiating both side, we have
\[
\frac{d}{dx}\left[ a(x)p(x)\right] = \frac{2b(x)}{a(x)}a(x)p(x).
\]
Hence, it holds that
\[
(\log a(x)p(x))' = \frac{2b(x)}{a(x)}.
\]
By integrate both sides from $c$ to $x$ we have the assertion.\qed

\vskip0.1cm

For $f \in C_0(S)$ (the set of continuous functions on $S$ vanishing at the boundary of $S$), let $m_f:= \int_{l}^{u} f(x)p(x)dx$ and define $\tilde g_f$ by, for every $x\in S$,
\begin{eqnarray}
\tilde g_f(x) &:=& \frac{2}{a(x)p(x)}\int _l ^x ( f(y)-m_f)p(y) dy \label{def-g1}\\
&=& - \frac{2}{a(x)p(x)}\int _x ^u ( f(y)-m_f)p(y) dy. \label{def-g2}
\end{eqnarray}
Then, by Proposition \ref{prop-ap} we have
\begin{equation}\label{eq-g}
\tilde g_f(x) = \int _l ^x \frac{2(f(y)-m_f)}{a(y)}\exp \left( -\int _y^x \frac{2b(z)}{a(z)}dz\right) dy, \hskip0.5cm x\in S.
\end{equation}
Then, $g_f(x):=\int _0^x \tilde g_f(y)dy$ satisfies that $f-m_f=Ag_f$ and by the definition of $m_f$ we have
\begin{equation}\label{eq2}
f(x) - \mathbf{E}[f(X)] = \frac12 a(x)\tilde g_f'(x) + b(x)\tilde g_f(x)
\end{equation}
where $X$ is a random variable with its law $\mu$.

\begin{remark}Relation (\ref{eq2}) is called Stein's equation and it is a general extension of the Stein's equation for normal and Gamma distributions considered in \cite{NoPe1}, \cite{NoPe4} or \cite{Reinert}. Indeed, when the measure $\mu$ is the standard normal distribution $N(0,1)$, the state space is $S=(-\infty, \infty)$ and the coefficients of the associated operator are defined by $a(x)= 2$ and $b(x)= -x$.
Therefore, (\ref{eq2}) becomes, for $X\sim N(0,1)$
\begin{equation*}
f(x)-\mathbf{E}f(X) = \tilde{g} _{f}'(x) -x \tilde{g}_{f}(x), \hskip0.5cm x \in \mathbb{R}.
\end{equation*}
When we deal with the Kolmogorov distance and take $f(x)= 1_{-(\infty, z] }(x)$ then
\begin{equation*}
\tilde{g}_{f}(x) =e^{\frac{x^{2}}{2}} \int_{-\infty}^{x} \left[ 1_{(-\infty, z]} (a) -\Phi (z) \right] e^{-\frac{a^{2}}{2}}da
\end{equation*}
which is the solution of the Stein's equation presented in \cite{NoPe1} or \cite{Reinert} for example.

We also retrieve the results in \cite{NoPe1} in the case when $\mu $ is the Gamma distribution. We refer to Section 4 for details.
\end{remark}

Now we consider the bounds for the functions  $\tilde g_f$ and $\tilde g_f'$.
\begin{prop}\label{prop-bound}
Assume that there exist $l', u' \in (l,u)$ such that $b$ is non-increasing on $(l,l')$ and $(u',u)$. Consider $f: S\to \mathbb{R}$ such that $\tilde{g}_{f}$ is well-defined and $\Vert f\Vert _{\infty}:= \sup_{x\in S} \vert f(x) \vert <\infty$.
Then we have
\begin{eqnarray*}
||\tilde g_f||_\infty \leq  C_1||f||_\infty  \mbox{ and }
||a\tilde g_f'||_\infty \leq C_2||f||_\infty
\end{eqnarray*}
 where $C_1$ and $C_2$ are strictly positive constants.
\end{prop}
{\bf Proof: }
Note that the condition imposed on $b$ implies $\liminf _{x\rightarrow u} b(x) < 0$ and $\limsup _{x\rightarrow l} b(x) > 0$.
By (\ref{cond-b}), (\ref{def-g2}) and (\ref{def-a}) we have
\[
|\tilde g_f(x)|= \left| \frac{\int _x^u(f(y)-m_f)p(y)dy}{\int _x^u b(y)p(y)dy} \right| .
\]
By L'H\^opital's rule we have
\[
\limsup _{x\rightarrow u}|\tilde g_f(x)|\leq  \limsup _{x\rightarrow u}\left| \frac{f(x)-m_f}{b(x)}\right| \leq C_1^+||f||_\infty
\]
where $C_1^+$ is a strictly positive constant.
Similarly we have
\[
\limsup _{x\rightarrow l}|\tilde g_f(x)| \leq C_1^-||f||_\infty
\]
where $C_1^-$ is a constant.
Hence, the continuity of $\tilde g_f$ yields the first assertion.
In view of (\ref{eq2}), to show the second assertion, it is sufficient to prove
\begin{equation}\label{proof-bound1}
||b\tilde g_f||_\infty \leq C_3 ||f||_{\infty}.
\end{equation}
By (\ref{cond-b}), (\ref{def-a}) and (\ref{def-g2}) we get
\[
|b(x)\tilde g_f(x)| \leq 2||f||_\infty \left| \frac{b(x)\int _x^up(y)dy}{\int _x^u b(y)p(y)dy}\right|.
\]
Since $b(x)$ is non-increasing on $[u',u)$, we have
\[
b(x)\int _x^up(y)dy \geq \int _x^u b(y)p(y)dy,\quad x\in [u',u).
\]
Hence, since $b$ is positive on $(k,u)$,
$
\limsup _{x\rightarrow u}|b(x)\tilde g_f(x)| \leq 2||f||_\infty .
$
and similarly we have
$
\limsup _{x\rightarrow l}|b(x)\tilde g_f(x)| \leq 2||f||_\infty .
$ Therefore, (\ref{proof-bound1}) holds from the continuity of $b\tilde g_f$. \qed

\begin{remark}
The hypotheses assumed on $b$ are satisfied for all the distributions considered throughout our paper. This is true because in all the examples the function $b$ is of the form $b(x)= -(x-m)$, $m$ being the expectation of the law $\mu$.
\end{remark}

The estimates in Proposition \ref{prop-bound} are sufficiently good when $a$ is uniformly bounded  and strictly positive.
But, when $a$ degenerates at the boundary of $S$, we need another estimate instead of the second estimate.
\begin{prop}\label{prop-bound2}
Assume that if $u<\infty$, there exists $u' \in (l,u)$ such that $b$ is non-decreasing and Lipschitz continuous on $[u',u)$ and $\liminf_{x\rightarrow u}a(x)/(u-x) >0$; if $u=\infty$, there exists $u' \in (l,u)$ such that $b$ is non-decreasing on $[u',u)$ and $\liminf_{x\rightarrow u}a(x)>0$.
Similarly, assume that if $l>-\infty$, there exists $l' \in (l,u)$ such that $b$ is non-increasing and Lipschitz continuous on $(l,l']$ and $\liminf_{x\rightarrow l}a(x)/(x-l) >0$; if $l=-\infty$, there exists $l' \in (l,u)$ such that $b$ is non-decreasing on $(l,l']$ and $\liminf_{x\rightarrow l}a(x)>0$.
Then we have
\[
||\tilde g_f'||_\infty \leq C_4 (||f||_\infty + ||f'||_\infty )
\]
for $f\in C_0^1(S)$ where $C_4$ is a constant.
\end{prop}
{\bf Proof: }
When $u=\infty$, we have $\limsup _{x\rightarrow \infty}|g_f'(x)| \leq C_2 ||f||_\infty$ by a similar argument to that in the proof of Proposition \ref{prop-bound}.
Similarly, when $l=-\infty$, we have $\limsup _{x\rightarrow -\infty}|g_f'(x)| \leq C_2 ||f||_\infty$.
In view of continuity of $g_f'$ it is sufficient to show that
\begin{eqnarray}
\limsup _{x\rightarrow u} |g_f'(x)| &\leq& C_5(||f||_\infty + ||f'||_\infty ), \quad \mbox{when}\ u<\infty , \label{proof-prop-bound2-u}\\
\limsup _{x\rightarrow l} |g_f'(x)| &\leq& C_5(||f||_\infty + ||f'||_\infty ), \quad \mbox{when}\ l>-\infty ,  \label{proof-prop-bound2-l}
\end{eqnarray}
where $C_5$ is a constant.
Consider the case that $u<\infty$.
Choose $\varepsilon >0$ such that $b(x)<-\varepsilon$ for $x\in [u',u)$ and $K>0$ such that $|b(x)-b(y)| \leq K|x-y|$ for $x,y\in [u',u)$.
By (\ref{cond-b}), (\ref{def-a}), (\ref{def-g2}) and (\ref{eq2}), we have
\[
\tilde g_f'(x) = \frac{2}{a(x)}\left\{ (f(x)-m_f)-\frac{b(x)\int _x^u (f(y)-m_f)p(y)dy}{\int _x^u b(y)p(y)dy}\right\} .
\]
Hence, for $x\in [u',u)$
{\allowdisplaybreaks
\begin{eqnarray*}
|\tilde g_f'(x)| &=& \frac{2}{a(x)\left| \int _x^u b(y)p(y)dy\right|}\left| (f(x)-m_f)\int _x^u b(y)p(y)dy -(f(x)-m_f)b(x)\int _x^u p(y)dy \right.\\
&& \phantom{\frac{2}{a(x)\left| \int _x^u b(y)p(y)dy\right|}} \ \left. + (f(x)-m_f)b(x)\int _x^u p(y)dy - b(x)\int _x^u (f(y)-m_f)p(y)dy \right| \\
&=& \frac{2}{a(x)\left| \int _x^u b(y)p(y)dy\right|} \left| (f(x)-m_f)\int _x^u [b(y)-b(x)]p(y)dy + b(x)\int _x^u [f(x)-f(y)]p(y)dy \right| \\
&\leq & \frac{2}{\varepsilon a(x)\int _x^u p(y)dy} \left[ 2K||f||_\infty (u-x) \int _x^u p(y)dy + b(x)||f'||_{\infty} (u-x)\int _x^u p(y)dy \right] \\
&=& \frac{2(u-x)}{\varepsilon a(x)} \left[ 2K||f||_\infty + b(x)||f'||_{\infty} \right] .
\end{eqnarray*}
}
This estimate and the assumption on $a$ implies (\ref{proof-prop-bound2-u}) with a constant $C_5$.
When $l>-\infty$, we can show (\ref{proof-prop-bound2-l}) by similar argument. \qed

There are many examples when our argument applies  in Table 1 on Page 8 in \cite{BSS}. Several examples will be discussed in details in the next section. We make below just some general comments.

We are now able to  derive the Stein's bound between the probability measure $\mu$ and the law of an arbitrary random variable $Y$. The following result extends the findings in \cite{NoPe1}, \cite{NoPe4} in the case where $\mu$ is standard normal law and the Gamma law respectively. We mention that $\langle \cdot, \cdot \rangle _{H}$ denotes the scalar product in ${\cal{H}}$.

\begin{theorem}\label{tt1}
Assume $X\sim \mu$ and let $Y$ be an $S$-valued random variable in $\mathbb{D}^{1,2}$ with $b(Y) \in L^2(\Omega)$. Then for every $f: S\to \mathbb{R}$ such that $\tilde{g}_{f}, \tilde{g}_{f}'$ are bounded,
\begin{eqnarray}|\mathbf{E}[f(Y) - f(X)]|
 &\leq &||\tilde g_f'||_{\infty} \mathbf{E}\left[ \left| \frac12 a(Y) + \langle D(-L)^{-1}\left\{ b(Y)-\mathbf{E}[b(Y)]\right\} , DY \rangle _H \right| \right] \nonumber\\
 &&+ ||\tilde g_f||_{\infty}|\mathbf{E}\left[ b(Y)\right] |.\label{eq3}
\end{eqnarray}
and
\begin{eqnarray}
|\mathbf{E}[f(Y) - f(X)]|
& \leq & ||\tilde g_f'||_{\infty} \mathbf{E}\left[ \left| \left. \mathbf{E}\left[ \frac12 a(Y) + \langle D(-L)^{-1}\left\{ b(Y)-\mathbf{E}[b(Y)]\right\} , DY \rangle _H \right| Y\right] \right| \right] \nonumber \\
&&+ ||\tilde g_f||_{\infty}|\mathbf{E}\left[ b(Y)\right] | . \label{eq3bis}
\end{eqnarray}
\end{theorem}
{\bf Proof: }
First,  by (\ref{eq2})
\begin{equation}\label{eq2.1}
\mathbf{E}[f(Y) - f(X)] = \mathbf{E}\left[ \frac12 a(Y)\tilde g_f'(Y) + b(Y)\tilde g_f(Y) \right] .
\end{equation}
 Recall that $LF:= -D^*DF$ for $F\in\mathbb{D} ^{1,2}$ centered. Then, since $(-L)^{-1}\{ b(Y)-\mathbf{E}[b(Y)]\} \in \mathbb{D} ^{1,2}$,
\begin{eqnarray*}
&& \mathbf{E}\left[ \frac12 a(Y)\tilde g_f'(Y) + b(Y)\tilde g_f(Y) \right] \\
&&= \mathbf{E}\left[ \frac12 a(Y)\tilde g_f'(Y) + \langle D(-L)^{-1}\left\{ b(Y)-\mathbf{E}[b(Y)]\right\} , \tilde g_f'(Y)DY \rangle _H \right] + \mathbf{E}\left[ b(Y)\right]\mathbf{E}\left[ \tilde g_f(Y)\right]
\end{eqnarray*}
Hence, by (\ref{eq2.1})
\begin{eqnarray}
&&\hspace{-1cm} |\mathbf{E}[f(Y) - f(X)]| \nonumber \\
&&\hspace{-1cm}  = \left| \mathbf{E}\left[ \tilde g_f'(Y) \left( \frac12 a(Y) + \langle D(-L)^{-1}\left\{ b(Y)-\mathbf{E}[b(Y)]\right\} , DY \rangle _H \right) \right] + \mathbf{E}[ \tilde g_f(Y)] \mathbf{E}\left[ b(Y)\right] \right| \label{eq3.1}\\
&&\hspace{-1cm}  \leq ||\tilde g_f'||_{\infty} \mathbf{E}\left[ \left| \frac12 a(Y) + \langle D(-L)^{-1}\left\{ b(Y)-\mathbf{E}[b(Y)]\right\} , DY \rangle _H \right| \right] + ||\tilde g_f||_{\infty}|\mathbf{E}\left[ b(Y)\right] |.\nonumber
\end{eqnarray}
It is possible to give an alternative bound for the distance.  Actually, the above calculation can be refined as follows.
\begin{eqnarray*}
&&\mathbf{E}[f(Y) - f(X)]\\
&=&\mathbf{E}\left[ \frac12 a(Y)\tilde g_f'(Y) + \langle D(-L)^{-1}\left\{ b(Y)-\mathbf{E}[b(Y)]\right\} , \tilde g_f'(Y)DY \rangle _H \right] + \mathbf{E}\left[ b(Y)\right]\mathbf{E}\left[ \tilde g_f(Y)\right] \\
&=&\mathbf{E} \left[ \tilde g_f'(Y)\mathbf{E}\left[ \left. \left( \frac12 a(Y) + \langle D(-L)^{-1}\left\{ b(Y)-\mathbf{E}[b(Y)]\right\} , DY \rangle _H \right) \right| Y\right] \right] + \mathbf{E}\left[ b(Y)\right]\mathbf{E}\left[ \tilde g_f(Y)\right]
\end{eqnarray*}and thus
\begin{eqnarray}
&&|\mathbf{E}[f(Y) - f(X)]| \nonumber \\
&& \leq ||\tilde g_f'||_{\infty} \mathbf{E}\left[ \left| \left. \mathbf{E}\left[ \frac12 a(Y) + \langle D(-L)^{-1}\left\{ b(Y)-\mathbf{E}[b(Y)]\right\} , DY \rangle _H \right| Y\right] \right| \right] + ||\tilde g_f||_{\infty}|\mathbf{E}\left[ b(Y)\right] | . \nonumber
\end{eqnarray}
\qed

\begin{remark}
a) Theorem \ref{tt1} is applicable for functions $f\in C_{0} ^{1}(S)$. Also the hypothesis $b(Y)\in L^2(\Omega)$ holds for $b(x)=-(x-m)$ once $Y\in \mathbb{D} ^{1,2}$.

b)  In the case $\mu =N(0,1)$ (recall $a(x)=2, b(x)=-x$) we obtain from Proposition 2 that for a centered random variable $Y\in \mathbb{D} ^{1,2}$
\[
\sup _{f\in C_0(S)} |\mathbf{E}[f(Y) - f(X)]| \leq C\mathbf{E}\left[ \left| 1 - \langle D(-L)^{-1} Y, DY \rangle _H \right| \right]
\]
where $C$ is a constant independent of $Y$. Notice that in this case we don't need $f\in C_{0}^{1}(S)$ due to Proposition 2 and Remark 4 since $\inf _{x\in \mathbb{R}} a(x)= 2$.  We mention also that is well known from \cite{Ste} that $\tilde g_{f}, \tilde{g}_{f}'$ are bounded under different assumptions on $f$ (see also Lemma 2.1 in \cite{NoPe1}).

\end{remark}

The Kolmogorov distance $d_K$ between the ${\mathcal L}(F)$ and ${\mathcal L}(G)$  (the laws of the random variables $F,G$) is defined by
\[
d_K({\mathcal L}(F), {\mathcal L}(G)) := \sup _{f\in {\mathcal F}_K} |\mathbf{E}[f(F)]-\mathbf{E}[f(G)]|.
\]
where ${\mathcal F}_K := \{ 1_{(l,z]}; z \in (l,u) \}$.
For $x\in S$ and $f(x)=1_{(l,z]}(x)$, we can choose $f_n \in C_0(S)$ such that $\{ f_n\}$ is an increasing sequence and $f_n(x)$ converges to $f(x)$ for all $x\in S$. Hence, by the dominated convergence theorem it holds that
\[
\lim _{n\rightarrow \infty} |\mathbf{E}[f_n(F)]-\mathbf{E}[f_n(G)]| = |\mathbf{E}[f(F)]-\mathbf{E}[f(G)]|.
\]
This implies the following estimate:
\[
d_K({\mathcal L}(F), {\mathcal L}(G)) \leq \sup _{f\in C_0(S);\ ||f||_\infty \leq 1}|\mathbf{E}[f(F)]-\mathbf{E}[f(G)]|.
\]
If $\inf _{x\in S}a(x)>0$, then by this estimate, Theorem \ref{tt1}, and Proposition \ref{prop-bound} we obtain an estimate for Kolmogorov distance between $X$ and $Y$ as follows:
\begin{eqnarray}
d_K({\mathcal L}(Y), \mu)   \leq C \mathbf{E}\left[ \left| \frac12 a(Y) + \langle D(-L)^{-1}\left\{ b(Y)-\mathbf{E}[b(Y)]\right\} , DY \rangle _H \right| \right]+ C|\mathbf{E}\left[ b(Y)\right] |,\label{eq100}
\end{eqnarray}
where $C$ is a positive constant.
Note that if $\mu$ is the normal distribution, we can choose $a(x)=2$ and $b(x)=-x$.

Generally, consider a distance between distributions of random variables $F$ and $G$ on $S$ defined by
\begin{equation}\label{general-distance}
d_{\mathcal H}({\mathcal L}(F), {\mathcal L}(G)) := \sup _{f\in {\mathcal H}} |\mathbf{E}[f(F)]-\mathbf{E}[f(G)]|,
\end{equation}
where ${\mathcal H}$ is a set of functions on $S$.
If for all $f\in {\mathcal H}$ there exists a sequence $f_n \in {\mathcal F}$ such that $f_n$ converges to $f$ in suitable sense, we have
\[
d_{\mathcal H}({\mathcal L}(F), {\mathcal L}(G)) \leq \sup _{f\in {\mathcal F}}|\mathbf{E}[f(F)]-\mathbf{E}[f(G)]|.
\]
Hence, by this estimate and (\ref{eq3}) we obtain an estimate for the distance between $X$ and $Y$ as follows:
\begin{eqnarray}&&d_{\mathcal H}({\mathcal L}(Y), \mu)
 \leq \sup _{f\in {\mathcal F}}||\tilde g_f'||_{\infty} \mathbf{E}\left[ \left| \frac12 a(Y) + \langle D(-L)^{-1}\left\{ b(Y)-\mathbf{E}[b(Y)]\right\} , DY \rangle _H \right| \right] \nonumber\\
 &&\hspace{5cm}+ \sup _{f\in {\mathcal F}}||\tilde g_f||_{\infty}|\mathbf{E}\left[ b(Y)\right] |.\label{eq101}
\end{eqnarray}

There are many kind of distance between distributions defined by (\ref{general-distance}).
For example,
by taking ${\mathcal H} = \{ f : ||f||_L \leq 1\}$, where $||\cdot ||_L$ denotes the usual Lipschitz seminorm, one obtains the Wasserstein (or Kantorovich-Wasserstein) distance;
by taking ${\mathcal H} = \{ f : ||f||_{BL} \leq 1\}$, where $||\cdot ||_{BL} = ||\cdot ||_L + ||\cdot ||_\infty$, one obtains the Fortet-Mourier (or bounded
Wasserstein) distance;
by taking ${\mathcal H}$ equal to the collection of all indicators $1_B$ of Borel sets, one obtains the total variation distance.

In the case of the Wasserstein distance or the Fortet-Mourier distance, for all $f\in {\mathcal H}$ we can choose ${\mathcal F}:= \{ f\in C_0^1(S); ||f'||_\infty \leq 1\}$, because we can choose $f_n \in \{ f\in C_0^1; ||f'||_\infty \leq 1\}$ such that $f_n$ converges to $f$ uniformly in every compact set.
Hence, (\ref{eq101}), Theorem \ref{tt1} and Proposition \ref{prop-bound2} implies that
\begin{eqnarray*}
d_{\mathcal H}({\mathcal L}(Y), \mu)  \leq C \mathbf{E}\left[ \left| \frac12 a(Y) + \langle D(-L)^{-1}\left\{ b(Y)-\mathbf{E}[b(Y)]\right\} , DY \rangle _H \right| \right]+ C|\mathbf{E}\left[ b(Y)\right] |,
\end{eqnarray*}
where $C$ is a positive constant.
Note that we do not assume $\inf _{x\in S}a(x)>0$ in this case.
This means that, even if $\mu$ has the Gamma distribution, the Wasserstein distance and the Fortet-Mourier distance are dominated.

In the case of the total variation distance, we choose ${\mathcal F}:= C_0(S)$, because we can choose $f_n \in C_0(S)$  uniformly bounded and such that $f_n(x)$ converges to $f(x)$ for each $x\in S$.
Hence, the similar argument to the case of the Kolmogorov distance is available, and if $\inf _{x\in S}a(x)>0$ we obtain the same estimate (\ref{eq100}) for the total variation distance.

We will discuss in the next section the significance of the bound given by (\ref{eq3}) and (\ref{eq3bis}).
The above computation leads to an interesting characterization of the random variations whose distribution is the invariant measure $\mu$ of the semigroup associated with the operator $A$.
\begin{theorem}
\label{char}
A random variable $Y \in \mathbb{D} ^{1,2}$ with its values on $S$ has probability distribution $\mu$ and satisfies that $b(Y)\in L^2(\Omega)$ if and only if $\mathbf{E}[b(Y)]=0$ and
\begin{equation}\label{eqchar}
\mathbf{E}\left[ \left. \frac12 a(Y) + \langle D(-L)^{-1} b(Y), DY \rangle _H \right| Y\right] =0.
\end{equation}
\end{theorem}
{\bf Proof: }
Suppose that $\mathbf{E}[b(Y)]=0$ and (\ref{eqchar}) holds. Then, due to (\ref{eq3}), the distance between the law of $Y$ and $\mu$ is zero and then $Y\sim \mu$. Suppose now that $Y\sim \mu$.
(\ref{cond-b}) implies that $\mathbf{E}[b(Y)]=0$.
Let $h \in C_K^\infty(S)$.
Define
\[
\tilde f(x):= \frac{1}{2p(x)}\frac{d}{dx}\left[ a(x)p(x)h(x)\right] ,\ x\in S .
\]
Note that (\ref{def-a}) implies that $ap\in C^1(S)$.
Since $h$ has compact support in $S$, we have $\tilde f\in C_0(S)$ and
$
m_{\tilde f} := \int _l^u \tilde f(x)p(x)dx =0.
$
Hence, the definition of $\tilde f$ and (\ref{def-g1}) implies that $h= \tilde g_{\tilde f}'$ where $\tilde g_{\tilde f}$ defined by $\tilde g_f$ with replacing $f$ by $\tilde f$.
This argument yields that for $h \in C_K^\infty(S)$ there exists $f\in C_0(S)$ such that $h= \tilde g_{\tilde f}'$.
Thus, (\ref{eq3.1}) implies
\[
\mathbf{E}\left[ h(Y) \mathbf{E}\left[ \left. \frac12 a(Y) + \langle D(-L)^{-1} b(Y) , DY \rangle _H \right| Y\right] \right] =0,\quad h\in C_K^\infty(S).
\]
This finishes the proof, because $C_K^\infty(S)$ is dense in $C_0(S)$ and the functions in $C_{0}(S)$ approximate the indicator functions. \qed

\begin{remark}
a) The same result has been obtained in \cite{NV} or \cite{T} in the case of the Gaussian distribution.

b) When the random variable $\langle DY, D(-L) ^{-1} (b(Y)-\mathbf{E}b(Y))\rangle _{H}$ is measurable with respect to the sigma algebra generated by $Y$ the bound (\ref{eq3}) and (\ref{eq3bis}) coincide. But, as it can be seen in the sequel, this is not always the case.
\end{remark}

\section{Significance of the bound}
The purpose of the section is to check the significance of the bound (\ref{eq3}). By "significance of the bound" we mean the following: given a random variable whose probability law is the invariant measure $X$, then the distance between its law and $X$ is zero. We will  prove that this is true in the case of several continuous probability distributions: the uniform distribution, the log-normal distribution and the Pareto distribution. But it fails in the case of the Laplace distribution. That means, for an explicit  random variable $Y$ which follows the Laplace distribution, we will prove that the right hand side of (\ref{eq3}) does not vanish almost surely.  On the other hand, as we have showed in Theorem \ref{char}, the right hand side of (\ref{eq3bis}) vanishes almost surely. This fact can be interpreted as follows: when the random variable $\langle DY, D(-L) ^{-1} b(Y) \rangle $ is not measurable with respect to the sigma-algebra generated by $Y$, the correct Stein's bound is the inequality (\ref{eq3bis}).

In order to compute the right hand side of (\ref{eq3}) we need to calculate the random variable $\langle DY, D(-L) ^{-1} b(Y) \rangle $. This random variable (and its conditional expectation given $Y$) appears in several works related to Malliavin calculus and Stein's method (see \cite{NoPe1}, \cite{T}, \cite{V}).  In general, it is difficult to find an explicit expression for it $Y$ for general $Y$. But in the case when $Y$ is a function of a Gaussian vector we have a very useful formula proved in \cite{NV}:
if  $Y= h(N)-\mathbf{E} h(N)$ where $h: \mathbb{R} ^{n}\to \mathbb{R}$   is a function of class $C^{1}$
with bounded derivatives and $N=(N_{1},..., N_{n}) $ is a Gaussian vector   with zero mean and covariance matrix $K=(K_{i,j}) _{i,j=1,..,n}$ then (we will omit in the sequel the index $H$ for the scalar product)
\begin{equation}
\langle D(-L) ^{-1} (Y-\mathbf{E}Y), DY\rangle = \int_{0} ^{\infty} e^{-u} du \mathbf{E}' \sum_{i,j=1} ^{n} K_{i,j} \frac{\partial h }{\partial x_{i}} (N) \frac{\partial h} {\partial x_{j}} (e^{-u}N+ \sqrt{1-e^{-2u}}N' ).\label{nv}
\end{equation}
Here  $N'$ denotes and independent copy of $N$ and we assume that $N,N'$ are defined on a product probability space
 $\left( \Omega \times \Omega ', {\cal{F}} \otimes {\cal{F}}, P\times P'\right)$  and $\mathbf{E}' $
 denotes the expectation with respect to the probability measure $P'$.
  Formula (\ref{nv}) is a consequence of the Mehler formula (see e.g. \cite{N}) and it has been proved in \cite{NV}, subsection 3.2.1.
 In the rest of this section the following context will prevail: $(W_{t}) _{t\in [0,T]}$ will denotes a standard Wiener process on $(\Omega , {\cal{F}}, P)$,  by $W(h)$
 we will denote the Wiener integral of $h\in L^{2}([0,T])$ with respect to  $W$ and $W'$ will be an independent Wiener process on a probability space $(\Omega ', {\cal{F}}', P')$.

 \subsection{The Gamma distribution}
The case of the Gamma distribution is already known. The Stein's bound (\ref{eq3}) has been obtained in \cite{NoPe1}, \cite{NoPe4} and already discussed in Section \ref{generaltheory}. We prefer to discuss it further in order to compare the bounds (\ref{eq3}) and (\ref{eq3bis}). We will consider the random variable
$$Y= W(h) ^{2}$$
which has Gamma distribution with parameters $a= \frac{1}{2}$ and $\lambda = \frac{1}{2}$.
This is actually the chi-square distribution and its associated coefficient are
$$ a(x)= 4x \mbox{ and } b(x)= -(x-1).  $$
 Proposition \ref{prop-bound2}  implies that (note the assumptions are satisfied by the functions $a,b$ above defined on the state space $(0, \infty)$
 \begin{equation}\label{bGamma}
 \left| \mathbf{E} f(X)-\mathbf{E}f(Y) \right| \leq C(||f||_\infty + ||f'||_\infty ) \mathbf{E}\left| 2Y-\langle DY, D(-L) ^{-1} (Y-1)\rangle \right| + \left| \mathbf{E}Y-1 \right|
 \end{equation}
 for any $f\in C_{0} ^{1}(0, \infty)$. Note that $\mathbf{E}Y=1$ for our choice of $Y$.

\begin{remark}
The bound (\ref{bGamma}) is a variant of inequalities (3.48), (3.49) in \cite{NoPe1} which are stated for the centered Gamma law and different classes of functions instead of $C_{0} ^{1}(0, \infty)$.
\end{remark}
 We can easily compute the scalar product $\langle DY, D(-L) ^{-1} (Y-1)\rangle $ using (\ref{nv}) with $h(x)= x^{2}$.
\begin{eqnarray*}
\langle DY, D(-L) ^{-1} (Y-1)\rangle&=& 2W(x) \int_{0} ^{1} da 2 \mathbf{E}' \left[ aW(h)+ \sqrt{1-a^{2} } W'(h) \right] \\
&=& 4W(h) ^{2} \int_{0}^{1} ada = 2 W(h) ^{2} = 2Y.
\end{eqnarray*}
We can notice that the random variable $\langle DY, D(-L) ^{-1} (Y-1)\rangle$ is measurable with respect to the sigma algebra generated by $Y$.
Therefore the bound (\ref{eq3}) and (\ref{eq3bis}) coincide and  Theorem \ref{char} provide an interesting characterization of the chi-square distribution random variables in terms of the Malliavin derivatives: that is, a random variable $Y$ has chi-square distribution
(with one degree of freedom) if and only if $\langle DY, D(-L) ^{-1} (Y-1)\rangle =2Y$ almost surely.

\subsection{The Uniform distribution}
We will discuss the case of the uniform distribution $U([a,b])$ with $\infty <a<b<\infty$. The density of this law is $p_{a,b}(x)= \frac{1}{b-a} 1_{[a,b]}(x)$ and  the mean of $p_{a, b}$ is $\frac{a+b}{2}$. We will actually restrict to the particular case $a=0, b=1$  and let
 $p(x):=p_{0,1}(x)= 1_{[0,1]}(x)$. This density  is associated  with the stochastic differential equation
\begin{equation*}
dX_{t}= -(X_{t} -\frac{1}{2}) dt + \sqrt{ X_{t}(1-X_{t}) } dW_{t} \end{equation*}
in the sense that the solution $X$ to the above equation is ergodic with invariant measure $\mu \sim U([0,1])$.
The diffusion coefficients $a$ and $b$ are in this case defined on $(0,1)$ and given by   (see Table 1 in \cite{BSS})
\begin{equation}
a(x)= x(1-x) \mbox{ and } b(x)=-(x-\frac{1}{2}). \label{ab-u}
\end{equation}
Let $Y$ be a random variable in the space $\mathbb{D} ^{1,2}$ such that
$$\mathbf{E}Y = \mathbf{E}U[0,1]= \frac{1}{2}.$$
In this case the Stein's bound  (\ref{eq3}) becomes, for any function $f\in C^{1}_{0}([0,1]) $ (we mention that $a$ satisfies the assumptions in Proposition \ref{prop-bound2} since $\lim _{x\to 1} \frac{a(x)}{1-x}=\lim _{x\to 0} \frac{a(x)}{x}=1$)
\begin{eqnarray*}
\left| \mathbf{E}f(Y) -\mathbf{E}f(U[0,1])\right| &\leq & C \left| \mathbf{E}\frac{1}{2} a(Y) - \langle D(-L) ^{-1} (Y-\frac{1}{2}), DY\rangle \right|\\
&=&C \left| \mathbf{E}\frac{1}{2}Y(1-Y) -\mathbf{E}\langle D(-L) ^{-1} (Y-\frac{1}{2}) , DY \rangle \right| .
\end{eqnarray*}
Let us check  "how good" is this bound on an example. Let $ f,g \in L^{2} ([0,T])$ such that
$$\Vert f \Vert _{L^{2}([0,T])}= \Vert g \Vert _{L^{2}([0,T])}=1
\mbox{ and } \langle f, g\rangle _{L^{2} ([0,T] )} =0.$$
Then $W(f)$ and $ W(g)$ are independent standard normal random variables.
Define the random variable $Y$ by
\begin{equation}
\label{Y1}
Y= e ^{ -\frac{1}{2} \left( W(f)^{2} + W(g) ^{2} \right)}.
\end{equation}
Then it is well-known that $Y$ has uniform distribution $U([0,1])$ since the random variable
 $ -\frac{1}{2} \left( W(f)^{2} + W(g) ^{2} \right)$ has exponential distribution with parameter 1. It is also clear that $Y\in \mathbb{D} ^{1,2}$.
Note also that $Y$ can be expressed as a function of the Gaussian vector $(W(f), W(g))$ whose covariance matrix is the identity matrix $I_{2}$.
Indeed, $Y=h(W(f), W(g))$ with
$$h: \mathbb{R} ^{2} \to \mathbb{R}, \hskip0.3cm h(x,y)=e^{-\frac{1}{2}(x^{2}+y^{2}) }.$$
The function $h$ satisfies the assumption in order to apply (\ref{nv}). Applying this formula to the random variable (\ref{Y1}) we get
\begin{eqnarray}
&&\langle D(-L) ^{-1} (Y-\frac{1}{2}), DY\rangle \nonumber \\
&&= W(f) e^{-\frac{1}{2} (W(f)^{2} + W(g)^{2})}\mathbf{E}'\int_{0} ^{\infty} e^{-u}du (e^{-u} W(f) + \sqrt{1-e^{-2u} }W'(f)  )\nonumber \\
&&\quad \times e^{-\frac{1}{2}\left[ (e^{-u} W(f) + \sqrt{1-e^{-2u}} W'(f)  ) ^{2} +(e^{-u} W(g) + \sqrt{1-e^{-2u}} W'(g)  ) ^{2} \right] } \nonumber  \\
&&\quad +W(g) e^{-\frac{1}{2} (W(f)^{2} + W(g)^{2})}\mathbf{E}'\int_{0} ^{\infty} e^{-u}du (e^{-u} W(g) + \sqrt{1-e^{-2u} }W'(g)  )\nonumber \\
&&\quad \times e^{-\frac{1}{2}\left[ (e^{-u} W(f) + \sqrt{1-e^{-2u}} W'(f)  ) ^{2} +(e^{-u} W(g) + \sqrt{1-e^{-2u} }W'(g)  ) ^{2} \right] }. \nonumber
\end{eqnarray}
We made the change of variable $e^{-u}=a$ and then
\begin{eqnarray}
\langle D(-L) ^{-1} (Y-\frac{1}{2}), DY\rangle
 &=&W(f) e^{-\frac{1}{2} (W(f)^{2} + W(g)^{2})}\mathbf{E}'\int_{0} ^{1} da (a W(f) + \sqrt{1-a^{2}} W'(f)  )\nonumber \\
&& \times e^{-\frac{1}{2}\left[ (a W(f) + \sqrt{1-a^{2}} W'(f)  ) ^{2} +(a W(g) + \sqrt{1-a^{2}} W'(g)  ) ^{2} \right] }\nonumber  \\
 && + W(g) e^{-\frac{1}{2} (W(f)^{2} + W(g)^{2})}\mathbf{E}'\int_{0} ^{1} da (aW(g) + \sqrt{1-a^{2}} W'(g)  )\nonumber \\
&& \times e^{-\frac{1}{2}\left[ (a W(f) + \sqrt{1-a^{2}} W'(f)  ) ^{2} +(a W(g) + \sqrt{1-a^{2}} W'(g)  ) ^{2} \right] }. \label{dd}
 \end{eqnarray}
 At this point we need the following lemma. It will be widely used throughout the paper.
\begin{lemma}\label{aux1}
Let $K\geq -1$, $C\in {\mathbb{R}}$ and $a\in (0,1)$. Suppose $Z\sim N(0,1)$. Then
\begin{equation*}
\mathbf{E}e^{ -K( C + \sqrt{1-a^{2}}Z  ) ^{2}} = \frac{1}{ \sqrt{1+2K (1-a^{2})}} e^{-C^{2}K \frac{1}{1+2K (1-a^{2}) }}
\end{equation*}
and
\begin{equation*}
\mathbf{E}\left( C + \sqrt{1-a^{2}} Z\right) e^{ -K( C + \sqrt{1-a^{2}}Z  ) ^{2}}=\frac{C}{\left( 1+ 2K (1-a^{2}) \right) ^{\frac{3}{2}}}e^{-C^{2}K \frac{1}{1+2K (1-a^{2}) }}.
\end{equation*}
\end{lemma}
{\bf Proof: }
The proofs for the two equations are almost same. We only show the second one.
\begin{eqnarray*}
&&\mathbf{E}\left( C + \sqrt{1-a^{2}} Z\right) e^{ -K( C + \sqrt{1-a^{2}}Z  ) ^{2}}\\
&&= \frac{1}{\sqrt{2\pi }} \int_{\mathbb{R}} (C+\sqrt{1-a^{2}}x)e^{-K(C+ \sqrt{1-a^{2}} x)^{2}}e^{-\frac{x^{2}}{2}}dx \\
&&= \frac{1}{\sqrt{2\pi }} e ^{-KC^{2}}e^{\frac{2K^{2}C^{2}(1-a^{2})}{1+2K(1-a^{2})}} \int_{\mathbb{R}}(C+\sqrt{1-a^{2}}x)e^{-\frac{1}{2}\left( 1+2K (1-a^{2})\right) \left( x+ \frac{2KC \sqrt{1-a^{2}}}{1+2K(1-a^{2})} \right)^{2}} dx \\
&&= \frac{1}{\sqrt{1+2K (1-a^{2})}}\left( C -\frac{2KC(1-a^{2})}{1+2K(1-a^{2})}\right)e^{-C^{2}K \frac{1}{1+2K (1-a^{2}) }}\\
&&= \frac{C}{\left( 1+ 2K (1-a^{2}) \right) ^{\frac{3}{2}}}e^{-C^{2}K \frac{1}{1+2K (1-a^{2}) }}.
\end{eqnarray*} \qed

\vskip0.1cm

Now, using Lemma \ref{aux1}, we have
\begin{eqnarray*}
&\mathbf{E}'e^{-\frac{1}{2}(a W(g) + \sqrt{1-a^{2}} W'(g)  )^{2}}=e^{-\frac{1}{2} \frac{a^{2}W(g) ^{2}}{(2-a^{2})  }} (2-a^{2})^{-\frac{1}{2}}&\\
&\mathbf{E}'(a W(f) + \sqrt{1-a^{2}} W'(f))e^{-\frac{1}{2}  (a W(f) + \sqrt{1-a^{2}} W'(f) ) ^{2}}=W(f)e^{-\frac{1}{2} \frac{a^{2}W(f) ^{2}}{(2-a^{2})  }} (2-a^{2})^{-\frac{3}{2}}a.&
\end{eqnarray*}
By inserting the above two identities in (\ref{dd}) we obtain
\begin{eqnarray}
\langle D(-L) ^{-1} (Y-\frac{1}{2}), DY\rangle \nonumber = (W(f)^{2}+ W(g) ^{2}) \int_{0} ^{1} e^{-\frac{1}{2} (W(f) ^{2}+ W(g) ^{2}) \frac{2}{2-a^{2}}} (2-a^{2} ) ^{-2}a da.
\end{eqnarray}
Since for any constant $c$
$$\frac{d}{da} e^{c \frac{1}{2-a^{2}} } =e^{c \frac{1}{2-a^{2}} } \frac{2ac} {(2-a^{2})^{2}}$$ we get (with $c= -(W(f) ^{2} + W(g) ^{2}) $)
\begin{eqnarray}
&&\langle D(-L) ^{-1} (Y-\frac{1}{2}), DY\rangle = \frac{1}{2}\left( e^{-\frac{1}{2} (W(f) ^{2}+ W(g) ^{2})}- e^{- (W(f) ^{2}+ W(g) ^{2})}\right).\label{gx}
\end{eqnarray}
On the other hand, using (\ref{ab-u}),
\begin{equation}\label{gx1}
a(Y)= Y(1-Y)
\end{equation}
and by (\ref{gx}) and (\ref{gx1}) we concluded that the right hand side of (\ref{eq3}) is zero.

\begin{remark}
It is interesting to note that $\frac{1}{2}a(Y)-\langle D(-L) ^{-1} (Y-\frac{1}{2}), DY\rangle $
 is zero  and  not only the  expectation of its absolute value. That is, this quantity is zero for every $\omega \in \Omega$.     We will also mention that in this case the     random variable $\langle D(-L) ^{-1} (Y-\frac{1}{2}), DY\rangle$
 is measurable with respect to the sigma-algebra generated by $Y$. Therefore the bounds (\ref{eq3}) and (\ref{eq3bis}) coincide.

\end{remark}

\subsection{The Beta distribution}
Using the computations in the previous paragraph, it is immediate to treat the case of the beta distribution with a particular choice of its parameters. Recall that the density of the beta distribution with parameters $\alpha, \beta  >0$ is
$
p_{\alpha , \beta } (x)= \frac{ \Gamma (\alpha + \beta) } {\Gamma (\alpha ) \Gamma (\beta) } x^{\alpha -1} (1-x) ^{\beta -1} 1_{(0,1) }(x).
$
The mean of this law is $\frac{\alpha }{\alpha + \beta}$ while the coefficient of the diffusion associated with the beta law are
$
a(x)= \frac{2}{\alpha + \beta } x(1-x) \mbox{ and } b(x)= -(x -\frac{ \alpha } {\alpha + \beta }).
$
Similarly to the case of the unform distribution we can check that $a$ and $b$ satisfies the hypothesis of Proposition (\ref{prop-bound2}).
We will restrict here to the special case $\alpha =\frac{1}{2}, \beta =1$. The distribution $\beta (\frac{1}{2} , 1)$  is a power-function distribution. It is well-known that if  $X\sim U[0,1]$ then $X^{2} \sim \beta (\frac{1}{2}, 1) $.
Consider the random variable
$Y= e^{- ( W(f) ^{2} + W(g) ^{2}) } . $
Then obviously $Y\sim \beta (\frac{1}{2}, 1) $. Using formula (\ref{nv}) with $h(x, y)= e^{-x^{2} + y^{2}}$, Lemma \ref{aux1} we get
\begin{equation*}
\langle DY, D(-L) ^{-1} (Y-\frac{1}{3}) \rangle = \frac{1}{2} a(Y).
\end{equation*}
Again the random variable $\langle DY, D(-L) ^{-1} (Y-\frac{1}{3}) \rangle $ is measurable with respect to $Y$.

\subsection{The log-normal distribution}\label{section-lognormal}
We analyze here the case of the lognormal distribution. Let us first review some basic properties of this probability  distribution. A random variable $Y$ has lognormal distribution
 with parameters $\delta $ and $\sigma ^{2}$ if $\log Y$ has normal distribution with mean $\delta $ and variance $\sigma ^{2}$.
  The density of the log-normal distribution with parameters $\delta $ and $\sigma ^{2}>0$ is
\begin{equation}\label{flognormal}
\frac{1}{\sqrt{2\pi \sigma ^{2}}x} e^{-\frac{1}{2\sigma ^{2}}(\log x -\delta )^{2}}1_{(0, \infty )} (x)
\end{equation}
and the coefficients of the associated diffusion are defined on $(0, \infty)$ and given by
$$b(x)= -(x-e^{\delta + \frac{1}{2}  \sigma ^{2}})$$
and
\begin{equation}
\label{aln}a(x)=\frac{2}{p(x)} \left( \Phi \left( \frac{ \log x-\delta }{\sigma }\right) - \Phi \left( \frac{ \log x-\delta }{\sigma }-\sigma \right)\right)
\end{equation}
(see \cite{BSS}, page 8)), where $\Phi $ denotes the cumulative distribution function of the standard normal law, $p$ is given by (\ref{flognormal}).
We check the significance of the bound in the case of the lognormal distribution with parameters $\delta =0$ and $\sigma=1$. The function $a$ satisfies $\lim _{x\to \infty} a(x)=\infty $ and $\frac{ a(x)}{x} = 2 e^{\frac{1}{2} (\log x)^{2}}\int_{\log x -1 }^{\log x} e^{-\frac{u^{2}}{2}}du \geq 2 $.

Let us consider the random variable
\begin{equation*}
Y= e^{W(h)}
\end{equation*}
where $h \in L^{2} ([0,T])$ has $L^{2}$ norm equal to 1. Then obviously $Y$ follows a lognormal law with mean $\delta =0$ and variance $\sigma=1$.
Let us first compute the scalar product
$$ \langle DY, D(-L) ^{-1} (b(Y) -Eb(Y)) \rangle =- \langle DY, D(-L) ^{-1} (Y-e^{\frac{1}{2} }) \rangle.$$
 Using formula (\ref{nv}) we get
 \begin{eqnarray*}
 \langle DY, D(-L) ^{-1} (Y-e^{\frac{1}{2} }) \rangle&=& \int_{0} ^{\infty} du e^{-u} e^{W(h)} \mathbf{E}' \left[ e^{e^{-u} W(h) + \sqrt{ 1-e^{-2u}} W'(h)}\right] \\
 &=& e^{W(h)}\int_{0} ^{1} da \mathbf{E}'\left[ e^{aW(h)  + \sqrt{ 1-a^{2}} W'(h) } \right]
\end{eqnarray*}
where $W'(h)$ denotes an independent copy of $W(h)$. Since
\begin{equation*}
\mathbf{E}'\left[ e^{aW(h)  + \sqrt{ 1-a^{2}} W'(h) }\right]= e^{aW(h)}e^{\frac{1}{2} (1-a^{2})},
\end{equation*}
we obtain
\begin{eqnarray*}
&& \langle DY, D(-L) ^{-1} (Y-e^{\frac{1}{2} }) \rangle =e^{W(h)} \int_{0} ^{1} da  e^{aW(h)} e^{\frac{1}{2} (1-a^{2})}\\
 &=&e^{W(h)} e^{\frac{1}{2}}e^{\frac{W(h)^{2}}{2}} \int_{0}^{1} da e^{-\frac{1}{2} (a-W(h) ) ^{2}}=
 e^{W(h)} e^{\frac{1}{2}}e^{\frac{W(h)^{2}}{2}}\int_{W(h)-1} ^{W(h)} e^{-\frac{x^{2}}{2}}dx
 \end{eqnarray*}
 where we used the change of variables $W(h)-a=x$.
Let us now compute
\begin{eqnarray*}
\frac{1}{2}a(Y)&= &e^{\frac{1}{2}} p(Y) ^{-1} \left( \Phi (\log (Y) ) -\Phi (\log (Y) -1) \right) \\
&=&  e^{\frac{1}{2}} Y e ^{\frac{\log (Y) ^{2}}{2}} \int_{\log(Y)-1} ^{\log (Y)} e^{-\frac{x^{2}}{2}} dx =  e^{\frac{1}{2}} e^{W(h)} e^{\frac{W(h)^{2}}{2}}\int_{W(h)-1} ^{W(h)} e^{-\frac{x^{2}}{2}}dx
\end{eqnarray*}
where we use the formula (\ref{flognormal}) for the density of the lognormal distribution and the expression of the operator $a$. We can see that
$$\frac{1}{2}a(Y)-  \langle DY, D(-L) ^{-1} (Y-e^{\frac{1}{2} })\rangle =0.$$
\begin{remark}
Again the quantity $\frac{1}{2}a(Y)-  \langle DY, D(-L) ^{-1} (Y-e^{\frac{1}{2} })\rangle $ vanishes and not only its expectation.
 Also $\langle DY, D(-L) ^{-1} (Y-e^{\frac{1}{2} })\rangle $  is measurable with respect to the sigma algebra generated by $Y$.
\end{remark}

\subsection{The Pareto distribution}
Let us recall some basic properties of the Pareto distribution  with parameter $\alpha >1$ (denoted in the following by $Pareto (\alpha )$).
 The probability density function of this law is
\begin{equation*}
p(x)= \alpha (1+x) ^{-\alpha -1}
\end{equation*}
and its   expectation  is $m=\frac{1}{\alpha -1}$. The functions $a$ and $b$ associated to the diffusion equation whose invariant measure is $Pareto (\alpha)$ are given by
\begin{equation*}
a(x)= \frac{2}{\alpha -1} x(1+x) \mbox{ and } b(x)= -(x-\frac{1}{\alpha -1}), \hskip0.5cm x\in (0, \infty).
\end{equation*}
It is standard to see that $a, b$ verifies the statement of Proposition \ref{prop-bound2}.

We recall a well-known fact: if the random variable $X$ follows a Pareto distribution with parameter $\alpha >1$ then
$$\log (Y+1) \sim Exp (\alpha).$$
Let us consider the same context as in the previous examples. That is, we are on a probability space $(\Omega, {\cal{F}}, P)$ and let $(W_{t})_{t\in [0,T]}$ be a Wiener process on this space.  Consider two orthonormal elements $h_{1}, h_{2} \in L^{2}([0,T])$. Then $W(h_{1}), W(h_{2})$ are independent standard normal random variables and
$$W(h_{1} )^{2}+ W(h_{2}) ^{2} \sim Exp (\frac{1}{2})= \Gamma (1, \frac{1}{2} ). $$
Consider the random variable
$$Y= e^{\frac{1}{4} (W(h_{1}) ^{2}+ W(h_{2}) ^{2})}-1.$$
Then, since $\frac{1}{4} (W(h_{1}) ^{2}+ W(h_{2}) ^{2})\sim \Gamma (1,2)= Exp (2)$, we can see that  $Y$ follows a Pareto distribution with parameter $\alpha =2$.  Clearly, we have
\begin{equation*}
\frac{1}{2} a(Y) = Y(1+Y) = e^{\frac{1}{2} (W(h_{1}) ^{2}+ W(h_{2}) ^{2})}-e^{\frac{1}{4} (W(h_{1}) ^{2}+ W(h_{2}) ^{2})}.
\end{equation*}
Using (\ref{nv}) with $h(x,y) = e^{\frac{1}{4} (x^{2} + y^{2})} -1$
we will get   (we recall that $W'$ is an independent copy of $W$, see the beginning of this section)
\begin{eqnarray*}
&&\langle DY, D(-L) ^{-1} (Y-1) \rangle \\
&=& \frac{1}{2} W(h_{1})  e^{\frac{1}{4} (W(h_{1}) ^{2}+ W(h_{2}) ^{2})} \int_{0}^{1} da \frac{1}{2} \\
&&
\mathbf{E}'\left[  (aW(h_{1})+ \sqrt{1-a^{2}}W'(h_{1})  )e^{\frac{1}{4} \left( (aW(h_{1})+ \sqrt{1-a^{2}}W'(h_{1})) ^{2}+  (aW(h_{2})+ \sqrt{1-a^{2}}W'(h_{2})) ^{2} \right)}\right]\\
&&+\frac{1}{2} W(h_{1})  e^{\frac{1}{4} (W(h_{1}) ^{2}+ W(h_{2}) ^{2})} \int_{0}^{1} da \frac{1}{2} \\
&&
\mathbf{E}'\left[  (aW(h_{2})+ \sqrt{1-a^{2}}W'(h_{2})  e^{\frac{1}{4} \left( (aW(h_{1})+ \sqrt{1-a^{2}}W'(h_{1})) ^{2}+  (aW(h_{2})+ \sqrt{1-a^{2}}W'(h_{2})) ^{2} \right)}\right]\\
\end{eqnarray*}
and by Lemma \ref{aux1} with $K= \frac{-1}{4}$ and $C= aW(h_{1}), aW(h_{2}) $ respectively,  we can write
\begin{eqnarray*}
\langle DY, D(-L) ^{-1} (Y-1) \rangle = Y(1+Y) = \frac{1}{2}a(Y).
\end{eqnarray*}

\subsection{The Laplace distribution: Failure of the bound (\ref{eq3})}
The Laplace distribution with parameter $\alpha >0$ (denoted by $Laplace (\alpha )$) is a continuous probability distribution with density $
p_{\alpha }(x)= \frac{\alpha}{2} e^{-\alpha \vert x\vert }, \mbox{ for every } x\in \mathbb{R}.$ The mean of the law is $m=0$ and the diffusion coefficients (\ref{difcoef}) are
\begin{equation}
\label{ab-laplace}
a(x)= \frac{2}{\alpha ^{2}} (1+ \alpha \vert x\vert ), \hskip0.3cm b(x)= -x.
\end{equation}
Is is known that if $X_{1}, X_{2}$ are two independent random variables such that $X_{1}\sim Exp (\alpha ), X_{2}\sim Exp (\alpha ) $ then $X_{1}-X_{2} \sim Laplace (\alpha )$.

Let us analyze the case of the Laplace distribution with parameter $\alpha =1$. In this case, (\ref{ab-laplace}) reduces to
\begin{equation*}
b(x)= -x \mbox{ and }  a(x)= 2 (1+|x| ).
\end{equation*}
Here the state space is whole real line $(-\infty, \infty)$ and we can apply Proposition \ref{prop-bound} in order to obtain the Stein's bound. Consider the random variable
$$Y = \frac{1}{2} \left( W(h_{1} ) ^{2} + W(h_{2} ) ^{2} - W(h_{3} ) ^{2}- W(h_{4} ) ^{2}\right)$$
where as above $h_{i}$ ($i=1,..,4$) are orthonormal functions in $L^{2}([0,T])$. Since
\begin{equation*}
 \frac{1}{2} \left( W(h_{1} ) ^{2} + W(h_{2} ) ^{2}\right) \sim Exp (1) \mbox{ and } \frac{1}{2}\left( W(h_{3} ) ^{2}+ W(h_{4} ) ^{2}\right) \sim Exp (1)
\end{equation*}
it can be easily seen that $Y\sim Laplace (1).$ It is easy to compute the quantity $\langle DY, D(-L)^{-1} b(Y)\rangle $ using formula (\ref{nv}). We obtain,
\begin{eqnarray}
\langle DY, D(-L)^{-1} b(Y)\rangle &=& -\langle DY, D(-L)^{-1} Y\rangle \nonumber\\
&=& -\frac {1}{2}  \left( W(h_{1} ) ^{2} + W(h_{2} ) ^{2} +W(h_{3} ) ^{2}+ W(h_{4} ) ^{2}\right).\label{pa1}
\end{eqnarray}
It is  obvious that is this case the difference $\frac{1}{2} a(Y)+\langle DY, D(-L)^{-1} b(Y)\rangle $  does not vanish almost surely. This signifies that the bound given by inequality (\ref{eq3}) is not good in this case and it has to replaced by (\ref{eq3bis}). Theorem \ref{char} ensures that the right hand side of (\ref{eq3bis}) vanishes almost surely. The reason why the bound (\ref{eq3}) fails is given by the fact that the random variable $\langle DY, D(-L)^{-1} b(Y)\rangle$ is not measurable with respect to the sigma-algebra generated by $Y$.

The case of the Laplace distribution can be discussed on an other examples. Consider four independent standard normal random variables $W(h_{1}), W(h_{2}), W(h_{3}), W(h_{4})$ and define
$$ Y_{1}= W(h_{1} ) W(h_{2} ) + W(h_{3} ) W(h_{4}). $$
Then again $Y_{1}$ follows a Laplace distribution with mean zero and variance 1.  And we can see that again  the expression (\ref{pa1}) holds for the random variable $Y_{1}$. Then
\begin{eqnarray*}
&&\frac{1}{2} a(Y_{1})+\langle DY_{1}, D(-L)^{-1} b(Y_{1})\rangle\\
 &&= 1+ \left|  W(h_{1} ) W(h_{2} ) + W(h_{3} ) W(h_{4})\right| - \frac {1}{2}  \left( W(h_{1} ) ^{2} + W(h_{2} ) ^{2} +W(h_{3} ) ^{2}+ W(h_{4} ) ^{2}\right)
\end{eqnarray*}
and this does not vanish. On the other hand, we know from Theorem \ref{char} that
$$\mathbf{E}\left( \left. \frac{1}{2} a(Y_{1})+\langle DY, D(-L)^{-1} b(Y_{1})\rangle \right| Y_{1} \right) =0.$$
in this way we will obtain some interesting (and somehow unexpected) identities for functions of the Brownian motion, which are difficult to be proven directly. That is
\begin{eqnarray*}
&&\mathbf{E} \left( \left. \frac {1}{2}  \left( W(h_{1} ) ^{2} + W(h_{2} ) ^{2} +W(h_{3} ) ^{2}+ W(h_{4} ) ^{2}\right) \right| W(h_{1} ) W(h_{2} ) + W(h_{3} ) W(h_{4})\right)\\
&& =1+ \left|  W(h_{1} ) W(h_{2} ) + W(h_{3} ) W(h_{4})\right|
\end{eqnarray*}
and
\begin{eqnarray*}
&&\mathbf{E} \left( \left. \frac {1}{2}  \left( W(h_{1} ) ^{2} + W(h_{2} ) ^{2} +W(h_{3} ) ^{2}+ W(h_{4} ) ^{2}\right) \right| \frac{1}{2} \left( W(h_{1} ) ^{2} + W(h_{2} ) ^{2} - W(h_{3} ) ^{2}- W(h_{4} ) ^{2}\right)\right)\\
&&=1+ \left| \frac{1}{2} \left( W(h_{1} ) ^{2} + W(h_{2} ) ^{2} - W(h_{3} ) ^{2}- W(h_{4} ) ^{2}\right) \right|.
\end{eqnarray*}

\section{Example}
We will illustrate the bound obtained via Stein's method through an example. Consider $(h_{i}) _{i\geq 0}$ a sequence of orthonormal elements of $L^{2}([0,T])$ and define for every $i\geq 1$
\begin{equation}
\label{xiLN}
X_{i} = e^{-\left( W(h_{i}) ^{2}-1\right) }
\end{equation}
(the minus sign is added in order to have finite expectation) and
\begin{equation}\label{ynln}
Y_{N}= \left( X_{1}....X_{N} \right) ^{\frac{1}{\sqrt{2N}}}= e^{-\frac{1}{\sqrt{2N}}\sum_{i=1} ^{N} (W(h_{i}) ^{2}-1)}.
\end{equation}
Then, applying the central limit theorem to $\frac{1}{\sqrt{2N}}\sum_{i=1} ^{N} (W(h_{i}) ^{2}-1)$, we have that $Y_{N}$ converges in distribution, as $N\to \infty$, to the log-normal distribution with mean zero and  variance equal to 1. Let us compute the bound given by the right hand side of (\ref{eq3}) and (\ref{eq3bis}).
Define $a$, $b$ and $p$ same as in Section \ref{section-lognormal}.
In this case we have
 \begin{equation*}
 \left| \mathbf{E} f(Y_{N} ) -\mathbf{E}f(X) \right|  \leq C  \mathbf{E}\left| \frac{1}{2} a(Y_{N})- \langle DY_{N}, D(-L) ^{-1} \left( Y_{N} -e^{\frac{1}{2}} \right) \rangle \right|+ C'   \left| \mathbf{E}(Y_{N}-e^{\frac{1}{2}}) \right| .
 \end{equation*}
 Since, with $Z\sim N(0,1)$
 $$\mathbf{E} Y_{N}= e^{\sqrt{\frac{N}{2}}} \left(  \mathbf{E} e ^{-\frac {1}{\sqrt{2N} }Z^{2}}\right) ^{N}=   e^{\sqrt{\frac{N}{2}}}  \left( \sqrt{ 1+ \frac{2}{\sqrt{2N}}}\right) ^{-N} $$
 we can see, by studying the asymptotic behavior as $N\to \infty $ of the above sequence, that
 $\sqrt{N} \mathbf{E} b(Y_{N}) \to _{N\to \infty} C $ with $C$ a strictly negative constant. We compute now $a(Y_{N})$ where $a$ is the function given by (\ref{aln}) with $\mu =e^{\frac{1}{2}}, \delta =0$ and $\sigma =1$. Denote by
\begin{equation*}
S_{N}= \sum_{i=1}^{N} W(h_{i} )^{2}\quad \mbox{and}\quad Z_{N}= \frac{1}{\sqrt{2N}} \sum_{i=1}^{N} (W(h_{i}) ^{2}-1).
\end{equation*}
We will have
\begin{eqnarray}
\frac{1}{2} a (Y_{N})&=& \frac{\mu  }{p(Y_{N})} \int_{\log Y_{N} -1}^{\log Y_{N}} \frac{1}{\sqrt{2\pi}} e^{-\frac{x^{2}}{2}} dx = e^{\frac{1}{2}} e ^{-Z_{N} }e^{\frac{1}{2} Z_{N} ^{2}}\int_{-Z_{N}-1}^{-Z_{N} }  e^{-\frac{x^{2}}{2}} dx \nonumber \\
&=& e^{\frac{1}{2} \left[Z_{N}-1\right] ^{2}} \int_{-Z_{N}-1}^{-Z_{N} }  e^{-\frac{x^{2}}{2}} dx=e^{\frac{1}{2} \left[Z_{N}-1\right] ^{2}} \int_{Z_{N}} ^{Z_{N}+1}  e^{-\frac{x^{2}}{2}} dx. \label{ayn}
\end{eqnarray}
Now, using  (\ref{nv}) with $h(x_{1},..., x_{N})= e^{-\frac{1}{\sqrt{2N}} \sum_{i=1}^{N} (x_{i}^{2}-1)}$
we can write
\begin{eqnarray}
&&\langle DY_{N}, D(-L) ^{-1} \left( b(Y_{N}) -\mathbf{E}b(Y_{N}) \right) \rangle = \langle DY_{N}, D(-L) ^{-1} \left( Y_{N} -e^{\frac{1}{2}} \right) \rangle \nonumber \\
&&= \sum_{i=1} ^{N} \sqrt{\frac{2}{N}} W(h_{i}) e^{-Z_{N}}\nonumber\\
&&\quad \times \int_{0}^{1} da  \sqrt{\frac{2}{N}} \mathbf{E}' \left( a W(h_{i})+ \sqrt{1-a^{2}} W'(h_{i})\right) e^{-\frac{1}{\sqrt{2N}} \left[ \sum_{i=1}^{N} (a W(h_{i})+ \sqrt{1-a^{2}} W'(h_{i}))^{2}-1\right] } \nonumber\\
&&=\sum_{i=1} ^{N} \sqrt{\frac{2}{N}} W(h_{i}) e^{-Z_{N}}\nonumber\\
&&\quad \times \int_{0}^{1} da  \sqrt{\frac{2}{N}}\mathbf{E}' \left( a W(h_{i})+ \sqrt{1-a^{2}} W'(h_{i})\right)e^{-\frac{1}{\sqrt{2N}} \left[  (a W(h_{i})+ \sqrt{1-a^{2}} W'(h_{i}))^{2}-1\right] } \nonumber\\
&&\quad \times \left( \prod_{j=1; j\not= i}^{N}\mathbf{E}' e^{-\frac{1}{\sqrt{2N}} \left[  (a W(h_{j})+ \sqrt{1-a^{2}} W'(h_{j}))^{2}-1\right] } \right).\label{dl1}
\end{eqnarray}
By applying Lemma \ref{aux1}  with $K=\frac{1}{\sqrt{2N}}$ and $C=aW(h_{i})$ we obtain
\begin{eqnarray*}
&&\mathbf{E}' e^{-\frac{1}{\sqrt{2N}} (a W(h_{j})+ \sqrt{1-a^{2}} W'(h_{j}))^{2} }= \frac{1}{ \sqrt{ 1+ \frac{2}{\sqrt{2N}}(1-a^{2})}}e^{-a^{2} W(h_{j}) ^{2}\frac{1}{\sqrt{2N}}\frac{1}{  1+ \frac{2}{\sqrt{2N}}(1-a^{2})}}\\
&&\mathbf{E}'\left( a W(h_{i})+ \sqrt{1-a^{2}} W'(h_{i})\right)e^{-\frac{1}{\sqrt{2N}}(a W(h_{i})+ \sqrt{1-a^{2}} W'(h_{i}))^{2}}\\
&&\hspace{5cm}=\frac{aW(h_{i})}{\left(  1+ \frac{2}{\sqrt{2N}}(1-a^{2})\right) ^{\frac{3}{2}}}e^{-a^{2} W(h_{i}) ^{2}\frac{1}{\sqrt{2N}}\frac{1}{  1+ \frac{2}{\sqrt{2N}}(1-a^{2})}}.
\end{eqnarray*}
Using the above two identities and (\ref{dl1})
\begin{eqnarray*}
&&\langle DY_{N}, D(-L) ^{-1} \left( Y_{N} -e^{\frac{1}{2}} \right) \rangle \\
&&=\frac{2}{N} \sum_{i=1}^{N} W(h_{i}) ^{2} e^{-Z_{N}} e^{ \sqrt{\frac{N}{2}}} \int_{0}^{1} da a \left( \prod_{i=1}^{N} e^{-a^{2} W(h_{i}) ^{2}\frac{1}{\sqrt{2N}}\frac{1}{  1+ \frac{2}{\sqrt{2N}}(1-a^{2})}}\right) \left( 1+ \sqrt{\frac{2}{N}}(1-a^{2})\right) ^{-\frac{N+2}{2}}\\
&&=\frac{2}{N} S_{N}  e^{-Z_{N}} e^{ \sqrt{\frac{N}{2}}} \int_{0}^{1} da a e^{ -a^{2} \frac{1}{\sqrt{2N}}S_{N}\frac{1}{  1+ \frac{2}{\sqrt{2N}}(1-a^{2})}}\left( 1+ \sqrt{\frac{2}{N}}(1-a^{2})\right) ^{-\frac{N+2}{2}}\\
&&=(\frac{2\sqrt{2} }{\sqrt{N}}Z_{N} +2) \int_{0} ^{1} da a e^{-Z_{N}(1+a^{2} \frac{1}{  1+ \frac{2}{\sqrt{2N}}(1-a^{2})})} e^{(1-\frac{a^{2}}{  1+ \frac{2}{\sqrt{2N}}(1-a^{2})})\sqrt{\frac{N}{2}}}\left( 1+ \sqrt{\frac{2}{N}}(1-a^{2})\right) ^{-\frac{N+2}{2}}.
\end{eqnarray*}
Since
\begin{equation*}
\log \left( 1 + \sqrt{\frac{2}{N}}(1-a^{2}) \right) = \sqrt{\frac{2}{N} } (1-a^{2})- \frac{1}{2} \frac{2}{N} (1-a^{2}) ^{2} +o(N^{-1})
\end{equation*}
where $o(N^{-1})$ is of Landau notation, and hence
\begin{equation*}
\left( 1+ \sqrt{\frac{2}{N}}(1-a^{2})\right) ^{-\frac{N+2}{2}}= e^{ -\frac{N+2}{2} (\sqrt{\frac{2}{N} } (1-a^{2})- \frac{1}{2} \frac{2}{N} (1-a^{2}) ^{2} +o(N^{-1})) },
\end{equation*}
we will have
\begin{eqnarray*}
&&\langle DY_{N}, D(-L) ^{-1} \left( Y_{N} -e^{\frac{1}{2}} \right) \rangle
=\frac{2\sqrt{2} }{\sqrt{N}}Z_{N} +2) \int_{0} ^{1} da a e^{-Z_{N}(1+a^{2} \frac{1}{  1+ \frac{2}{\sqrt{2N}}(1-a^{2})})} \\
&& \times e^{(1-\frac{a^{2}}{  1+ \frac{2}{\sqrt{2N}}(1-a^{2})})\sqrt{\frac{N}{2}}}e^{ -\frac{N+2}{2} (\sqrt{\frac{2}{N} } (1-a^{2})- \frac{1}{2} \frac{2}{N} (1-a^{2}) ^{2} +o(N^{-1})) }= A_{N} + B_{N},
\end{eqnarray*}
where
\begin{eqnarray*}
A_{N}&:=&\frac{2\sqrt{2} }{\sqrt{N}}Z_{N} \int_{0} ^{1} da a e^{-Z_{N}(1+a^{2} \frac{1}{  1+ \frac{2}{\sqrt{2N}}(1-a^{2})})} e^{(1-\frac{a^{2}}{  1+ \frac{2}{\sqrt{2N}}(1-a^{2})})\sqrt{\frac{N}{2}}}e^{ -\frac{N+2}{2} (\sqrt{\frac{2}{N} } (1-a^{2})- \frac{1}{2} \frac{2}{N} (1-a^{2}) ^{2} +o(N^{-1})) }\\
B_{N}&:=&2 \int_{0} ^{1} da a e^{-Z_{N}(1+a^{2} \frac{1}{  1+ \frac{2}{\sqrt{2N}}(1-a^{2})})} e^{(1-\frac{a^{2}}{  1+ \frac{2}{\sqrt{2N}}(1-a^{2})})\sqrt{\frac{N}{2}}}e^{ -\frac{N+2}{2} (\sqrt{\frac{2}{N} } (1-a^{2})- \frac{1}{2} \frac{2}{N} (1-a^{2}) ^{2} +o(N^{-1})) }.
\end{eqnarray*}

We will first show that
\begin{equation*}
\mathbf{E}\left|\sqrt{N} A_{N} \right| \to _{N\to \infty} C_{0}
\end{equation*}
where $C_{0} $ is a strictly positive constant.  We can write, using the change of variables $a^{2} =b$ with $2ada=db$
\begin{eqnarray*}
\mathbf{E}\left|\sqrt{N} A_{N} \right|&=& \sqrt{2} \mathbf{E}|Z_{N} | \int_{0} ^{1} db e^{-Z_{N}(1+ b \frac{1}{ 1+ \sqrt{\frac{2}{N}}(1-b)})} e^{\sqrt{\frac{N}{2}}\left[ 1 -b \frac{1}{ 1+ \sqrt{\frac{2}{N}}(1-b)}-(1-b) \right] } e^{g_{N}(b) }
\end{eqnarray*}
where $g_{N}(b)= e^{\frac{1}{2} (1-b) ^{2} + o(1)}.$ Therefore
\begin{eqnarray*}
\mathbf{E}\left|\sqrt{N} A_{N} \right|&=&\sqrt{2} \mathbf{E}| Z_{N} | \int_{0}^{1} db e^{-Z_{N}(1+b) } e^{b(1-b)}e^{g_{N} (b)} \\
&&+ \sqrt{2}\left(\mathbf{E}|Z_{N} | \int_{0} ^{1} db e^{-Z_{N}(1+ b \frac{1}{ 1+ \sqrt{\frac{2}{N}}(1-b)})} e^{\sqrt{\frac{N}{2}}\left[ 1 -b \frac{1}{ 1+ \sqrt{\frac{2}{N}}(1-b)}-(1-b) \right] } e^{g_{N}(b) } \right. \\
 &&\left. -    \mathbf{E}| Z_{N} | \int_{0}^{1} db e^{-Z_{N}(1+b) } e^{b(1-b)} e^{g_{N} (b)}\right)
\end{eqnarray*}
Since $Z_{N}$ converges in law to $N(0,1)$, the first summand above converges as $N\to \infty$ to
\begin{equation*}
C_{0}=\sqrt{2} e^{-1} \mathbf{E}\left| Z\right| \int_{0}^{1} db e^{-Z(1+b) }e^{b(1-b)}e^{\frac{1}{2} (1-b) ^{2}}
\end{equation*}
while the second summand converges to zero as $N\to \infty$ by the dominated convergence theorem.
Let us handle the summand denoted by $B_{N}$. Actually, we will prove that
\begin{equation}\label{bn-a}
\sqrt N \mathbf{E}\left| B_{N} -\frac{1}{2} a(Y_{N}) \right| \to _{N\to \infty} D_{0}
\end{equation}
with $D_{0}$ a strictly positive constant. First, we write
\begin{eqnarray*}
B_{N}&=&\int_{0}^{1}db e^{-Z_{N}(1+ b \frac{1}{ 1+ \sqrt{\frac{2}{N}}(1-b)})}e^{\sqrt{\frac{N}{2}}\left[ 1 -b \frac{1}{ 1+ \sqrt{\frac{2}{N}}(1-b)}-(1-b) \right] } e^{g_{N}(b) }\\
&=& \int_{0}^{1} db e^{-Z_{N} (1+b)} e^{b(1-b)} e^{\frac{1}{2} (1-b) ^{2}}\\
&&+ \int_{0}^{1} db \left( e^{-Z_{N}(1+ b \frac{1}{ 1+ \sqrt{\frac{2}{N}}(1-b)})}e^{\sqrt{\frac{N}{2}}\left[ 1 -b \frac{1}{ 1+ \sqrt{\frac{2}{N}}(1-b)}-(1-b) \right] } e^{g_{N}(b) } -e^{-Z_{N} (1+b)} e^{b(1-b)} e^{\frac{1}{2} (1-b) ^{2}}\right)\\
&:=& B_{N} ^{(1)}+ B_{N} ^{(2)}
\end{eqnarray*}
and the limit (\ref{bn-a}) follows since (see (\ref{ayn})) $
B_{N} ^{(1)} = \frac{1}{2} a(Y_{N})$ and for $N$ large,
$$1 -b \frac{1}{ 1+ \sqrt{\frac{2}{N}}(1-b)}-(1-b)= b(1-b) \sqrt{\frac{2}{N}}+ o(\frac{1}{\sqrt N})$$
and
$$ \left| e^{-Z_{N}(1+ b \frac{1}{ 1+ \sqrt{\frac{2}{N}}(1-b)})} -e^{-Z_{N} (1+b)} \right| \leq c \frac{|Z_N|}{\sqrt N}.$$
As a conclusion of the computations contained in this section, the distance between the law of $Y_{N}$ given (\ref{ynln}) by and the log-normal distribution with mean 0 and variance 1 is of order of $\frac{1}{\sqrt{N}}$.


\begin{thebibliography}{99}
\bibitem{Barbour}
{A.D. Barbour (1990): }\emph{ Stein's method for diffusion approximations. } Probability Theory and Related Fields 84, 297-322.
\bibitem{AS}
{Y. A\"it-Sahalia (1996): }{\em Nonparametric pricing of interest  rate derivative securities. } Econometrica, 64, 527-560.
\bibitem{BS}
{B.M. Bibby and M. Sorensen (2001): }{\em Simplified estimating functions for diffusion models with a high-dimensional parameter. } Scandinavian Journal of Statistics, 28(1), 99-112.
\bibitem{BSS}
{B.M. Bibby, I.M. Skovgaard and M. Sorensen (2003): }{\em Diffusion-type models with given marginals and auto-correlation function. } Bernoulli, 11(2), 191-220.

\bibitem{EV}
{R. Eden and F. Viens (2010): }\emph{General upper and lower tail estimates using Malliavin calculus and Stein's equations. } Preprint.

\bibitem{KaTa}
{S. Karlin and H. Taylor (1981): } A second course on stochastic processes. Academic press.
\bibitem{NoPe1}
{I. Nourdin and G. Peccati (2007): }\emph{ Stein's method on Wiener
chaos.} Probability Theory and Related Fields \textbf{145} (1-2), 75-118.
\bibitem{NoPe2}
{I. Nourdin and G. Peccati (2009): } {\emph Stein's method and exact Berry-Ess\'{e}en asymptotics for functionals of Gaussian fields. } The Annals of Probability,  37(6), 2200-2230.
\bibitem{NoPe3}
{I. Nourdin and G. Peccati (2008): }Stein's method meets Malliavin calculus: a short survey with new estimates.  \emph{Recent Advances in Stochastic Dynamics and Stochastic Analysis, } World Scientific.
\bibitem{NoPe4}
{I. Nourdin and G. Peccati (2009): } \emph{Non-central convergence of multiple integrals,} Ann. Probab. 37, no. 4, 1412-1426.
\bibitem{NV}
{I. Nourdin and F. Viens (2009): }{\em Density formula and concentration inequalities with Malliavin calculus. } Electronic Journal of Probability, 14, paper 78,  2287-2309.
\bibitem{N} {D. Nualart (2006): }\emph{Malliavin Calculus and Related Topics. Second Edition. }{Springer. }
\bibitem{Reinert}
{G. Reinert } {\em A short introduction to Stein's method. } Lecture Notes.
\bibitem{Ste}
{Ch. Stein (1986): } \emph{Approximate computation of expectations. } Institute of Mathematical Statistics Lecture Notes-Monograph Series \textbf{7}, Hayward, CA.

\bibitem{T} {C.A. Tudor (2011): } {\em Asymptotic Cram\'er theorem and analysis on Wiener space. } Seminaire de Probabilit\'es, Lectures Notes in Mathematics.
    \bibitem{V} {F. Viens (2009): } {\em Stein's lemma, Malliavin calculus, and tail bounds, with application to polymer fluctuation exponent. } Stoc. Proc. and their Applic. 119, 3671-3698.
\end{thebibliography}
\end{document}